\documentclass[10pt]{article}

\usepackage[dvips]{graphicx}
\usepackage{amsmath,amsfonts,amssymb}
\usepackage{a4wide}

\newtheorem{Proposition}{Proposition}[section]
\newtheorem{Theorem}{Theorem}[section]
\newtheorem{Definition}{Definition}[section]
\newtheorem{Remark}{Remark}[section]
\newtheorem{Corollary}{Corollary}[section]
\newtheorem{Lemma}{Lemma}[section]
\newtheorem{Example}{Example}[section]

\def\N{{\mathbb N}}

\def \dim {{\rm dim\,}}

\def \deg {{\rm deg\,}}
\def \E {{\cal E}}

\def \Ra {\Rightarrow}
\def \P {\mathbb P}

\def \M {${\cal M}_\Gamma$}

\def \Ra {\Rightarrow}

\def \H {{\mathfrak H}}
\def \I {{\mathfrak I}}
\def \M {{\mathfrak M}}
\def \J {{\mathfrak J}}
\def \A {{\mathfrak A}}
\def \E {{\mathfrak E}}
\def \L {{\mathfrak L}}
\def \p {{\mathfrak P}}
\def \basis {${\cal B}(\I)$}
\def \ffdd {\vrule height 4pt width 3pt depth 2pt}
\def \finedim {\par \nobreak \rightline{\ffdd \qquad}}

\begin{document}

\title{ \textbf{Essential and inessential elements of a standard basis}}

\author{\\ Giannina Beccari, Carla Massaza \\
Dipartimento di Matematica,  Politecnico Torino ,Italy\\}

\maketitle 

\date{}

In this paper we introduce the concept of {\it inessential element} of a standard basis \basis,  where $\I$ is any homogeneous ideal of a polynomial ring. An inessential element is, roughly speaking, a form of \basis \ whose omission produces an ideal having the same saturation of $\I$; it becomes useless in any dehomogenization of $\I$ with respect to a linear form. We study the properties of \basis\ linked to the presence of inessential elements and give some examples.

\section{Introduction}

The systems of generators of a given ideal $\I \subset K[y_0,...,y_n] =S$, satisfying given conditions, are widely studied. In the special case of homogeneous ideals, it is well known that there exist systems of generators, called standard bases, satisfying the following condition: their elements of degree $d$ are forms defining a $K$-basis of the vector space $\I_d/(\I_{d-1}S_1)$, for every $d \in \N$ (\cite{DGM},\cite{Ca},\cite{ Ro}). The standard bases are minimal among the systems of generators of a homogeneous ideal, but they are not the only interesting ones (for instance, Gr\"{o}bner   bases are not, in general, minimal, but they are of interest for other reasons). However, in this paper we will consider only standard bases of homogeneous ideals. The elements of each of them may be of two different kinds: essential generators (e.g.) and inessential generators (i.g.). A generator $g$ is called inessential, with respect to a basis \basis containing it, if it lies in the saturation of the ideal generated by ${\cal B}(\I)-\{g\}$;  this means that any dehomogenization $\I_*$ of $\I$ with respect to a linear form  is generated by the image of ${\cal B}(\I)-\{g\}$. A generator $g$ not lying in the saturation of ${\cal B}(\I)-\{g\}$ is called essential. We needed this concept in our attempt of considering the elements of a standard basis of $\I$ as separators (\cite{BM}) with respect to a convenient ideal $\J\subset \I$;  in fact, we found that such an interpretation is possible iff the generators are essential. We realized that not all the ideals do have a basis whose elements are all essential, neither in the special situation  of ideals of 0-dimensional schemes, in which we were interested. This fact suggested to study the concept of essentiality, independently from its use in the link between separating sequences and generators of a sub-ideal and for ideals of any height. So, we pointed our attention not only on the standard bases with the maximum number of essential elements, that at first had interested us, but also on those whose inessential elements are all contained in the saturation of the ideal generated by the essential ones:  any dehomogenization of $\I$ is generated by the images of their essential elements.
\medskip

Section 2 contains recalls and notation.

\smallskip
Section 3 contains the definition of essential and inessential elements, with equivalent formulations and some examples.

\smallskip
In section 4 we consider the special case of perfect height 2 ideals. In this situation, the essentiality or the inessentiality of a generator can be read as a property of the ideal generated by the entries of its corresponding column in any Hilbert matrix of $\I$ (\cite{N}).

\smallskip

In section 5 we come back to the study of the general situation. We show that any saturated homogeneous ideal has at least a basis with the maximum (resp. minimum) number of essential generators in any degree; so, the two sequences of those numbers are numerical sequences linked to the ideal; their elements are, degree by degree, less than or equal to the corresponding graded Betti numbers. We will call those bases e-maximal (resp. e-minimal) and give an algorithm of construction of one of them starting from any standard basis. An e-maximal (resp. e-minimal) basis is characterized by the fact that its inessential (resp. essential) elements have their typical property with respect to every standard basis containing them. The e-maximal bases were the first object of our interest, as we were looking for bases with the greatest number of generators to be viewed as elements of a separating sequence. We give just a few examples of search of e-maximal bases, as we are planning to devote to them another paper, in which we study a family of perfect height 2 ideals, for which it is possible to compute the number of the essential elements contained in an e-maximal basis, starting from some properties of their generators in minimal degree. From another point of view, a minimal e-basis seems to be of interest when we dehomogenize with respect to a linear form; in fact, an inessential element becomes useless as a generator of the dehomogenized ideal. However, from this point of view it turns out to be more suitable the notion of  "inessential set", generalizing the one of inessential element. In fact, the standard bases giving rise to a basis of minimal cardinality, after a dehomogenization with respect to a generic linear form, are the ones containing an inessential set of maximal cardinality. So, the last part of section 5 is devoted to such bases and to the ones (E-bases) whose set of inessential elements is an inessential set.

\section{Recalls and Notation}

Let $S=K[y_0,...,y_n]$, $K$ algebraically closed, be the coordinate ring of $\P^n$, $\I =\bigoplus \I_d, d\in \N$, a homogeneous ideal of $S$, $\M = (y_o,...,y_n)$ the irrelevant ideal. We recall the following:

\begin{Definition} \label {Def. 0.1}{\rm  (\cite{DGM}) \  A  {\it standard basis} \basis \ of $\I$ is an ordered set of forms of $S$, generating $\I$, such that its elements of degree $d$ define a $K$-basis of $\I_d/(\I_{d-1}S_1)$.}

\end{Definition}

It is well known(\cite {DGM}) that the number of generators of \basis \ , in a given degree $d$, depends only on the ideal $\I$: it is the $d$-th Betti number of $\I$, at the first level.

When we need to point out a subset $T$ of \basis \ , we use the non-standard notation:

${\cal B}(\I) = (t_1,...,t_m, s_1,...,s_p)$, where $T=(t_1,...,t_m)$ and $S=(s_1,...,s_p)$ inherit the ordering of \basis,\ which is, however, considered as an ordered set, with its original ordering.

\smallskip

Moreover,when there is no matter of misunderstanding, we will use the notation $(f_1,...,f_r)$ to denote the ideal generated by the standard basis $(f_1,...,f_r)$, instead of the heavier notation $(f_1,...,f_r)S$.

\medskip
If $\I$ is perfect of height $2$ (\cite{Ca}, \cite{N}, \cite {DGM}) , it is useful to consider, for every basis \basis, \ a {\it Hilbert  matrix} (\cite{N}), as follows.

Let: \begin{description}
\item ${\cal B}(\I) =(g_1,...,g_t)$, deg $g_j \leq $deg $g_{j+1}$,
\item $s_i = ( a_{i1},...,a_{it})$  $i$th element of a basis of syzygies with respect to \basis, \ where deg $s_i$=deg($a_{ij}g_j), j=1...t$,\  deg $s_i< $deg$s_{i+1}$,
 \end{description}

then $M(\I) = (a_{ij}), i=1...t-1,j=1...t$ \ is a Hilbert matrix of $\I$, related to \basis. \ Moreover (\cite {N}) $(-1)^j g_j$  is the minor of $M(\I)$ obtained by deleting its $i$-th column $C_j$. We say that $g_j$ is the generator linked to the column $C_j$ or that $C_j$ is its corresponding column.

The ideal generated by $C_j$ will be denote $\I_{C_j}$.

 \medskip

 We will be mainly interested in saturated ideals. We recall that:

 \begin{Definition}\label{Def. 02} {\rm (\cite{ Ha}) A homogeneous ideal $\I$ is saturated iff:
 $$ (\exists t) f \M^t \in \I  \Rightarrow  f\in \I.$$
 Equivalently, we can say that the irrelevant ideal is not associated to $\I$.}
 \end{Definition}

 Every ideal $\I$ has a saturation $\I^{sat}$, which is the minimum saturated ideal containing it. A process of computation of $\I^{sat}$, starting from $\I$, can be found in (\cite{Ha}).

 To every projective scheme $V$ of $\P^n$ we can associate a unique saturated ideal $\I$, which is usually denoted $\I(V)$.

  For every linear form $L\in S$,  the ideal $\I_*$, obtained from a homogeneous ideal $\I$ by dehomogenization with respect to $L$, is the image of $\I$ in the localization of S with respect to $L$. With a change of coordinate, it is possible to choose  $L=y_0$; in this situation, the localization of $S$ is isomorphic to $R = K[x_1,...,x_n]$ under the map associating to every  form $F(y_0,...,y_n)\in S$ the polynomial  $F_*(x_1,...,x_n) = F(1,x_1,...,x_n)\in R$ (see (\cite{F}) and $I_*$ can be identified with the image of $\I$ under that morphism .Viceversa, the homogenization  $\J^* \subset S$ of any ideal $\J\subset R$ is the ideal generated by $F(x_1,...,x_n)^* = y_0^d F(y_1/y_0,...,y_n/y_0)$, where $F(x_1,...,x_n)$ is any polynomial of $\I$ and $d$ is its degree. Let us observe that $(\J^*)_*=J$, while $(\I_*)^*= \I$ only if $y_0$ is regular for $S/\I$.
  The operation of dehomogenization can be made on a set of generators of $\I$, but the analogous is not true for the homogenization.

  Let us recall the following:

  \begin{Definition}\label{Def.03}{\rm A basis of an affine ideal $\J$ is a set of generators that fails to generate $\J$ if one of its elements is omitted.}
  \end{Definition}

 Two different bases of $\J$ may have a different cardinality.

\section{ Equivalent conditions and examples}

Let $\I$ , $\I^{sat}\varsubsetneqq M$,   be a homogeneous ideal of $S = K[y_0,...,y_n]$, $f$ an
element of a standard basis \basis \ and  $\H _{(f,{\cal B} )}$ the ideal generated
by  ${\cal B}(\I) - \{f \}$.

\medskip
\begin{Definition} \label {Def. 2.2} {\rm An element $f$ of a standard
basis ${\cal B}(\I)$ of  $\I$ is called  {\it inessential
generator of} $\I$ {\it with respect to} ${\cal B}(\I)$ iff  $f
\in (\H_{(f,{\cal B})})^{sat}$.

Otherwise, we say that $f$ is an {\it essential generator of }$\I$
{\rm with respect to } ${\cal B}(\I)$.}

\end{Definition}

\medskip
The following proposition gives conditions equivalent to inessentiality.

\begin{Proposition}\label {Prop. 2.1} {\rm  Let $\H = (g_1,...,g_r)$ be a homogeneous ideal of $S =
K[y_0,...,y_n]$, such that  $\H^{sat }\varsubsetneqq M =(y_0,...,y_n)$, $f$ any form of $\M$, $\I  =
(\H, f)$.   The following facts are equivalent:}
\begin {description}
 \item i) {\rm There exists $t\in \N$ such that $ f\M^t \subseteq \H$
(in other words, $ f\in \H^{sat}$).}

 \item ii) {\rm $S/\I$ and $S/\H$ have the same Hilbert
 polynomial (see \cite{Ha},\cite{Ro}).}

 \item iii) {\rm There exists a linear form $z\in S$, regular for
 $S/\H^{sat}$, such that a dehomogenization with respect to $z$
 gives: $\H_* = \I_*$.}

 \item iv)  {\rm $\H_* = \I_*$, for every dehomogenization with respect
 to any linear form $z \in \M $.}

\end {description}

\end {Proposition}

\medskip \noindent
{\bf Proof}

 The equivalence between $i)$ and $ii)$ is obvious.

$i) \Ra iv)$  \  Condition $i)$ implies $fz^t \in \H$, for every
$z\in \M$; as a consequence, $f_* \in \H_*$ or, equivalently,
$\H_* = \I_*$.

$iv \Ra iii)$  \ The condition $\H^{sat} \not= \M$ assures the
existence of an element $z$ regular for $S/\H^{sat}$ (as the union of its
associated primes cannot be $\M$),  so that the implication is obvious.

$iii) \Ra i)$ \ $f_* \in \H_*$ means $f_* = \Sigma \alpha_i g_{i*}
$, so that $(\exists t) fz^t = \Sigma \beta_i g_i \in \H
\subset \H^{sat}$.  As $z$ is regular for $S/ \H^{sat}$, we get \ $f\in
\H^{sat}$.

\finedim

\bigskip
\noindent {\bf Remarks}

1. \ Let us observe that it is sufficient to verify condition
$iv)$ for a set of linear forms generating $\M$; equivalently, the condition $f\in \H^{sat}$ is verified iff, for every linear form $L$ of a set of generators of $\M$, there exists $n\in \N$ such that $fL^n \in \H$.
\smallskip

2. \ In $iii)$ the condition \lq \lq $z$ is regular for $S/\H^{sat}$\rq \rq \ cannot
be replaced by  \lq \lq $z$ is regular for $S/\I$\rq \rq , as we can see in the
following  example.

Let:  $S = K[x,y,z], \  \I = (g_1, g_2,g_3,g_4) , \ \H =
(g_1,g_2,g_4), \ f=g_3$, \ where:

$$ g_1 = x^5,\quad g_2 = xy^5, \quad g_3 = y^7, \quad g_4 =
y^3(-x^4-y^2z^2)$$
 are the maximal minors of the matrix:

$$\begin{pmatrix} 0  &  x^3  &  z^2  &   -y^2 \\
                  0  &  y^2  &  -x   &     0  \\
                  y^3&  z^2  &   0   &    -x   \end{pmatrix}.$$

\smallskip

It is easy to verify that $z$ is regular for $S/\I$ and that, in
the dehomogenization with respect to $z$, we have: \ $ g_{3* }=
y^7 = -y^2 g_{4*}-x^3 g_{2*}$, that implies : $ g_{3*} \in \H_*$.
\ However, $g_3$ does not satisfy the equivalent conditions of
Proposition \ref{Prop. 2.1} ( see Proposition \ref {Prop. 2.9} for a quicker check). \
The reason is that $z$ is not regular for $\H^{sat}$. In fact: \
$z^2y^5 = -g_4 -x^4y^3$, \ so that $z(zy^7) = -g_4y^2-g_2x^3  \in
\H \subseteq \H^{sat}$  \ and \ $zy^7 \notin \H^{sat}$, \ as $zy^k \notin
\H, \forall k\in \N$.

\smallskip
3. \ To compute $\H^{sat}$ may be uneasy, so that condition $iv)$ and
Remark 1. become of some interest.

\finedim



\bigskip

 We give a statement equivalent to the essentiality of $f\in \I_d$, with respect to a standard basis ${\cal B}(\I) = \{f,{\cal B}_1\}$.
 \medskip
\begin{Proposition} \label{Prop. 2.17}{\rm Let  ${\cal B}=\{f, {\cal
B}_1\}$ be a standard basis of an  ideal $\I \subset
S=K[y_0,...,y_n], \ f\in \I_d$. The following facts
are equivalent:

$i)$ $f$ is essential with respect to ${\cal B}$;

$ii)$ there exists a set $\{L^*_1,...,L^*_n,N\}$ of linear forms,
generating $\M$, such that:

\qquad a) $N$ is a regular form both for $S/\H^{sat}, \ \H = ({\cal B}_1)
S$ and for $S/\I^{sat}$.

\qquad b) $\forall t \in \N, f N^t \notin\J,$ where $\J = ({\cal
B}_1, L^*_1 f,...,L^*_n f)$ .}
\end{Proposition}

\medskip
\noindent {\bf Proof}

\qquad $i)\Ra ii)$ \ As the union of the primes associated to $\H^{sat}$ or to $\I^{sat}$ cannot coincide with $\M$, we choose $N\in \M$ regular both for $S/\H^{sat}$ and $S/\I^{sat}$. Let us remark that:

(+)\quad $N$ does not divide $f$; otherwise $ f = f_1 N^u \in \I$,
$N$ regular for $S/\I$, would imply $f_1 \in \I$, so that $f$
would not satisfy the condition $ f \notin (\I_{d-1})S$.

A dehomogenization with respect to $N$ gives that $\I_* = (
f_*,\H_*)$ and $\H*$ are different ( see Proposition \ref{Prop. 2.1}
, iii)) and, as a consequence, there exists a maximal ideal
$\p=(L_1,...,L_n) \in R = S_*$, such that:
$$(*)\qquad \I_* R_{\p} \not= \H_* R_{\p}.$$
 Let us set $\J=(\H, fL_1^*,...,fL^*_n)$  and consider  $\J_* =
(\H_*, f_*L_1,...,f_*L_n)$; we prove that:
$$\J_* R_{\p} \not= \I_* R_{\p}.$$

Otherwise, in $R_{\p}$, we should have:
$$ f_* = (\sum_{i=1}^n a_i L_i) f_* +h_*,\ h_* \in \H_* R_{\p},$$
   so that $(1-\sum_{i=1}^n a_i L_i) f_*  \in \H_* R_{\p}$. \ As
$1-\sum_{i=1}^n a_i L_i$ is invertible in $R_{\p}$, this
implies $f_* \in \H_* R_{\p}$, against $(*)$. As a consequence, we
have $\I_* \not= \J_*$. Hence, \ $(\forall t) \ f N^t \notin \J$,\
for, otherwise, we should find, by dehomogenization with respect
to $N$: \ $f_* \in \J_*$,\ which implies $\I_* = \J_*.$

\qquad $ii) \Ra i)$ \ This implication comes immediately, as b)
implies:  $f N^t \notin ({\cal B}_1) S, \ \forall t$  and, as a
consequence,  $ (\forall t)  \ f \M^t \not\subset \H$; so, condition $ii)$ b)\ is sufficient to imply $i)$.

\finedim

\medskip
\noindent  {\bf Remark 1} \qquad Clearly it is enough to verify
condition \  b) for $t\gg 0$.

\noindent {\bf Remark 2}\qquad   If  $\dim S/\I =1$, condition b) can be
replaced by the following:

\qquad b')  $f$ is a {\it separator} for $S/\J$.( see[  ]).

In fact, in this case  the definition of
 separator is meaningful and condition b) can be restated as:

$\dim_K ( R/\J)_t = \dim_K(R/\I)_t +1, \  t \geq d.$

This relation, with condition (+), is equivalent to say that $f$
is a separator for $S/\J$.
\finedim

\medskip

With the same notation of Proposition \ref{Prop. 2.17}, we can state:

\medskip

\begin{Proposition} \label{Prop. 2.17'}{\rm If $\I$ is saturated, then so is $\J$.}
\end {Proposition}

\medskip
\noindent {\bf Proof}
\smallskip
Let us prove that $\J = \J^{sat}$. If not, we could find an element $u\in\M, u\notin J$ and a number $s\in \N$ such that $ u\M ^s\subset\J \subset\I$. As $\I$ is saturated, $u$ must be in $\I$; as a consequence, $u= aN^tf + j, a \in K^*, j \in \J$. But this implies $fN^t\M^s \subset \J \subset \I$, a contradiction, as $N$ is regular for $S/\I$ and $f$ is essential.
\finedim

\bigskip
The essentiality of an element $f$ depends on the basis in which
it is considered, as we can see in the following

\begin{Example} \label{Ex. 2.3}\quad {\rm Let $\I \subseteq K[x,y,z]$ be
the ideal generated by the maximal minors  of the matrix:

$$ M = \begin{pmatrix} z  &  0  &  0  &   -x \\
                      0  &  x  &  0  &    -y \\
                      0& 0  &  y   &    -z  \end{pmatrix}.$$

A standard basis of $\I$ is ${\cal B} (\I) = (g_1,g_2,g_3,f)$,
where:
$$ g_1 = x^2y, \ g_2= y^2z, \  g_3 = xz^2, \ f = xyz.$$

 We can easily
check that $f$ is inessential with respect to ${\cal B}(\I)$. In
fact: $fx = x^2yz = zg_1, \ fy = xy^2z = xg_2, \ fz = xyz^2 =
yg_3$,  \ so that  $f \M \subset (g_1, g_2, g_3) = \H$.

Let us produce a new basis, with respect to which $f$ is
essential. Choose in $\P^2$ a point not lying on $f = 0$, for
instance $P(1,1,1)$, and replace $g_1, g_2, g_3$ with generators
vanishing at $P$, so obtaining the new basis  ${\cal B'}(\I) =
(g_1-f, g_2-f, g_3-f, f)$. Clearly $\H' =( g_1-f, g_2-f, g_3-f)$
is such that $ (\H') ^{sat} \not = \I$, as the underlying schemes
differ for one point.

We observe that it is also possible to produce a standard basis
${\cal B}"(\I)$ such that every element of it is inessential: it
is enough to replace, in $M$, the first three columns with the sum
of each of them with the fourth (see Proposition \ref {Prop.
2.9}).}
\end {Example}

\finedim

\bigskip
There are also situations in which a generator of $\I$ is
inessential with respect to any basis containing it and every
basis contains at least an inessential element in the degree of
$f$. We see that situation in the following

\medskip
\begin{Example}\label{Ex. 2.4} \quad  {\rm Let $\I \subset K[x,y,z]$ be
the (saturated) ideal generated by the maximal minor of the
matrix:

$$M = \begin{pmatrix} y^2 &  0  &  0  &   -x \\
                      0  &  z^2  &  0  &    -y \\
                      0& 0  &  x^2   &    -z  \end{pmatrix}.$$

We have: \  $ {\cal B} (\I) = (g_1, g_2, g_3, f)$, \ where:

$$ g_1 = x^3z^2, \quad g_2 = x^2y^3, \quad g_3 = y^2z^3, \quad  f =
x^2y^2z^2.$$

It is immediate to see that \ $fx = y^2g_1,\ fy = z^2 g_2, \ fz =
x^2g_3$,\quad so that  \ $f\M \subseteq\H = (g_1, g_2, g_3)$.

In this case, every standard basis of $\I$ must contain an element
$f'$ in degree $6$, giving rise to the same ideal $\H$; moreover,
$f'$ must satisfy the relation $f' = kf + h,\quad k\in K, h\in\H$.
As a consequence, $f'\M \subseteq \H$ is still verified, so that
$f'$ is inessential with respect to any basis containing it.}
\end{Example}
\finedim

\noindent Taking into account the situation described in Example
\ref{Ex. 2.4}, we give the following

\medskip
\begin{Definition} \label{Def. 2.5} {\rm An element $f \in \I_d$ is {\it strongly inessential}
 ({\it s.i.}) iff $f\notin (\I_{d-1})S$  and it is inessential with respect to any standard
 basis containing it.}
\end{Definition}

\medskip
\begin{Proposition} \label{Prop. 2.6} {\rm   A strongly inessential
generator cannot have the minimal degree $\alpha (\I)$. }
\end{Proposition}

\medskip \noindent
{\bf Proof}
\smallskip
   Let us consider a basis \basis $=( f, h_1,...,h_r), \quad \deg \ f
   \leq \deg \ h_1 \leq ...\leq \deg \ h_r$. We prove that, if $f$ is
   inessential, it is possible to replace each $h_i$ with an
   $h'_i$, such that the behaviour of $h'_i$ with respect to
   essentiality is equal to the one of $h_i$ and $f$ is essential
   with respect to  ${\cal B}'(I) = (f,h'_1,...,h'_r)$.  To this
   aim, we choose a point $P\in {\mathbb P}^n$ and a linear form $z$ such
   that: \ $f(P) \not= 0,\ z(P) \not= 0$. \ In every linear
   system $h_i + \lambda _i z^{t_i} f, \ t_i = \deg f_i - \deg f, \ \lambda_i \in K$
    \ there is a form $h'_i = h_i + a_i z^{t_i} f $,
    \ such that \ $h'_i(P) = 0$. As a consequence, $f$ is essential
    with respect to ${\cal B}'(\I)$.
\finedim

\bigskip
\begin{Corollary} \label{Cor. 2.7} {\rm If $\I$ is generated in
minimal degree, then $\I$ admits a standard
basis of essential elements.}
\end{Corollary}

{\bf Proof}
\smallskip
 If ${\cal B}(f_1,...,f_r)$ is any standard basis and if $f_i$ is its first inessential generator, then, thanks to Proposition \ref{Prop. 2.6}, we can find a basis ${\cal B}'(f'_1,...,f'_i,...f'_r),  f'_j= f_j+a_j f_i, a_i=0$, with respect to which $f_i=f'_i$ is essential. Moreover,it is easy to check that$f'_j$ is still essential for $j<i$ ( see  Lemma \ref{Lemma. 2.1}). So, at any step, the basis ${\cal B}$ can be replaced by another basis with one more  essential element.

\finedim

\medskip

With a reasoning very similar to the one of Proposition \ref{Prop.
2.6}, we can prove:

\begin{Proposition} \label{Prop. 2.8} {\rm Let \basis
$=(h_1,...,h_m,f,g_1,...,g_k), \deg \ h_1 \leq...\leq \deg \ h_m <
\deg \ f \leq \deg \ g_1 \leq ... \leq \deg \ g_k$, $f$
inessential with respect to \basis.  If there exists a point $P
\in \P^n$ such that $h_i(P) = 0,\ i=1,...,m, \ f(P) \not= 0$, then
there exist $g'_1,...,g'_k$ such that $f$ is essential for ${\cal
B}'(\I) = (h_1,...,h_m, f, g'_1,...,g'_m)$. }
\end{Proposition}
\smallskip
\noindent Let us observe that the requirement of Proposition
\ref{Prop. 2.8} implies $ f\notin (h_1,...,h_m)^{sat}$; on the other
side, $f\in (h_1,...,h_m)^{sat}$ implies $f$ $s.i.$, but the viceversa
is not true, as we will see in Example \ref{Ex. 2.12}, where $g_3
\notin (g_1, g_2)^{sat}$, and $g_3$ is $s.i.$.

\bigskip

\noindent More generally, we would like to face the following problems:

 \begin{description}

\item {A.} Given a standard basis ${\cal B}(\I)$, find all its
elements of a given degree which are essential with respect to it.

\item{B.} Check how the "nature" (essentiality-inessentiality) of $f$ varies with the basis containing it.

\item {C.} Check how the number of essential elements in a given
degree varies with the chosen basis .

\end{description}

\section{The case of perfect height 2 ideals.}

\bigskip
If $\I$ is a perfect codimension $2$ ideal (for instance, the
ideal of a 0-dimensional scheme in $\P^2$), we can give an answer
to both questions A. and B. in terms of a Hilbert-Burch matrix
$M(\I)$ with respect to \basis.  If $f_r$ is the $r$-th element of
\basis, let us denote $\I_{C_r}\subset S$ the ideal generated by the
entries of the r-th column of $M(\I)$. With this notation, we can
state:

\medskip
\begin{Proposition} \label{Prop. 2.9}{\rm Let $\I$ be a perfect
codimension $2$ ideal of $S$. Then $f_r\in {\cal B}(\I)$ is
inessential for  \basis \ iff the following condition is
satisfied:

\begin{equation}
 \exists t \in \N, \ \M^t\subseteq \I_{C_r}  \label{eq. 2.1}
\end{equation}
}
\end{Proposition}

\medskip \noindent
{\bf Proof}

\smallskip
From the definition of $\I_{C_r}$ we get:
\begin{equation}
\I_{C_r} f_r \subseteq \H  \label{eq. 2.2}
\end{equation}
Conditions (\ref{eq. 2.1}) and (\ref{eq. 2.2}) imply:
\begin{equation}
\M^t f_r \subseteq \H,  \label{eq. 2.3}
\end{equation}
which says that $f_r$ is inessential.

Viceversa, (\ref{eq. 2.3}) implies the existence of syzygies whose
r-th components generate $\M^t$, so that $\M^t \subseteq \I_{C_r}$.
\finedim

\smallskip

\begin{Remark} \label{Re. 2.10} {\rm Proposition \ref{Prop. 2.9} can be restated
reducing the problem to the affine situation. Let $L$ be any
linear form of $S$, regular for $S/\H^{sat}$ and let $\I_*, \H_*, \I_{C_r*}$ be
the dehomogenization of $\I, \H, \I_{C_r}$. Then $f_r$ is inessential iff
$\I_{C_r *} = R$.}
\end{Remark}

\bigskip
Let us pass to consider problem B.. It is well known that a change
of a standard basis \basis $=(g_1,...,g_r,...g_m)$ is equivalent
to a change of its matrix $M(\I)$, realized by repeatedly
replacing a column $C_r$ with $ C'_r = \sum_{i} t_{ir}C_i$, where
$T = (t_{ir})$ is an invertible matrix and $ t_{rr} \in K^*,\
t_{ir} =0$ if $\deg \ g_i  < \deg \ g_r$. So, Proposition
\ref{Prop. 2.9} gives rise to the following

\medskip
\begin{Corollary} \label{Cor. 2.11}
{\rm Let $g_r \in {\cal B}(\I)$ be the generator corresponding to
the column $C_r$ of $M(\I)$; $g_r$ is $ s.i$ iff the entries of
every $C'_r = \sum_{j} t_{rj} C_j,  \ t_{rr} \in K^*$, generate an ideal
$\I_{C'_r}$ satisfying condition (\ref{eq. 2.1}) of Proposition
\ref{Prop. 2.9}}.

\end{Corollary}

\medskip \noindent
Let us consider again Examples  \ref{Ex. 2.3}  and \ref{Ex. 2.4}
from this point of view.
\smallskip
In Example \ref{Ex. 2.3}, the ideals generated by the entries of
its columns $C_i, i=1,...,4$ are respectively:

$\I_{C_1}= (z), \ \I_{C_2} = (x), \ \I_{C_3} = (y), \I_{C_4} = (x,y,z)$.

The only one satisfying the condition of Proposition \ref{Prop.
 2.9} is $\\I_{C_4}$, so that the only inessential element is $f$.

 Now, let us replace $C_4$ with a new column $C'_4$, so that the
 fourth generator becomes essential. We have :

 $$^t C'_4 = (-x+t_1z \qquad -y+t_2x \qquad  -z+t_3y)\qquad \ t_i \in K.$$

The ideal $\I_{C'_4}$ generated by $ C'_4$'s entries cannot contain a
power of $\M$ iff the linear system:

$$\begin{cases} -x+t_1z = 0 \\
              -y+t_2x = 0 \\
              -z+t_3y = 0
\end{cases}$$
has proper solutions, that is iff $t_1 t_2 t_3 = 1$. In
particular, choosing  \ $t_1 = t_2 = t_3 = 1$, we find again the
basis  ${\cal B}'(\I)$ already obtained with another technique.

\smallskip
In Example \ref{Ex. 2.4} the column corresponding to $f$ is the
fourth; it cannot be changed (apart from the multiplication by a
scalar) by degree reason: so, we find again that any standard
basis has an inessential generator in degree $6$.

\smallskip
Let us observe that, in this example, the inessential generator is
the only generator of maximal degree, so that its corresponding
ideal $\H$ does not depend on the standard basis. In the following
example the considered ideal $\I$ has, in every standard basis, an
inessential element of degree $11$, even if $11$ is not the
greatest degree of its generators, and another inessential element
in the maximal degree, in which there are two generators.
\medskip
\begin{Example} \label{Ex. 2.12}{\rm Let $\I \subset K[x,y,z]$ be
the ideal generated by the maximal minors of the matrix:

$$M(\I) = \begin{pmatrix} 0 &  x^5  &  0    &   -y^3  &   0 \\
                         0  &   0   &  x^3  &    z^2  &  -y^2 \\
                         0  &   0   &  y^2  &    -x   &   0  \\
                         y^3&   0    &  z^2   &    0   &   -x  \end{pmatrix}.$$

We have \ \basis $= (g_i), i=1...5$, where:

$$ g_1 = x^{10}, \ g_2 = y^{10}, \ g_3 = x^6 y^5, \  g_4 = x^5y^7, \
g_5 = x^5y^5z^2 + x^9y^3 .$$

The ideals generated by the entries of the columns $C_i, \
i=1...5,$ are respectively:

$$ \I_{C_1} = (y^3), \ \I_{C_2} = (x^5), \ \I_{C_3} = (x^3,y^2,z^2), \  \I_{C_4} =
(y^3, z^2, x), \ \I_{C_5} = (x, y^2).$$

So: $\I_{C_3}\supset\M^5, \ \I_{C_4}\supset \M^4$,  while $\I_{C_1}, \I_{C_2},
\I_{C_5}$ do not contain any power of $\M$. As a consequence, $g_1,
g_2, g_5$ are essential, while $g_3$ and $g_4$ are inessential.
After a general basis change of $\I$   modifying only the third
column, $C_3$ is replaced by $C'_3 = k C_3 + P C_4 + Q C_5$, where
$k\in K^*$, $P$ and $Q$ are linear forms. So:

$$\I_{C'_3} = (P y^3, kx^3 + P z^2 - Q y^2,  k y^2 - P x,  k z^2 - x
Q).$$
 It is immediate to control that the generators of $\I_{C'_3} $
are annihilated only by $ x = y = z = 0$, for every choice of  $k,
P, Q$. This means that the only prime ideal associated to $\I_{C'_3}$
is $\M$, so that $g_3$ is still inessential. Analogously, after a
general change of basis modifying only the fourth column, $C_4$ is
replaced by $C'_4 = k_1 C_4 + k_2 C_5, \ k_1 \in K^*, k_2 \in K$,
so that :
$$\I_{C'_4} = (-k_1 y^3, k_1 z^2 - k_2 y^2, -k_1 x, -k_2 x) = (x,
y^3, k_1z^2- k_2 y^2).$$
 Also in this case the generators of $\I_{C'_4}$ are annihilated only
 by $x=y=z=0$, so that $g_4$ is still inessential.}
\end{Example}

\medskip
Finally, let us give an example of an ideal with two $s.i.$\
generators in maximal degree.

\smallskip
\begin{Example} \label{Ex. 2.13} {\rm  Let $\I \subset K[x,y,z]$
be the ideal generated by the maximal minors of the matrix:

$$M(\I) = \begin{pmatrix} z^2 &  0  &  0    &   y  &   x \\
                         0  &  x^2+y^2   &  0  &    0  &  y \\
                         0  &   0   &  0  &    x   &  z  \\
                         0  &   0   &  x^2-y^2   &   z   &  0  \end{pmatrix}.$$

Any linear combination of the last two columns produces a new
column generating $\M$; this means that the two generators of
maximal degree are strongly inessential. More precisely, $\I =
(g_i)i=1...5$, where: $g_1=(x^4-y^4)(x^2-yz),\
g_2=xyz^2(x^2-y^2),\ g_3=z^4(x^2+y^2),\ g_4=z^3(x^4-y^4), \
g_5=xz^2(x^4-y^4)$ \
 and $\L = (g_1, g_2, g_3)$ is such
that $\L^{sat} = \I$, as can be seen with a direct computation; however, the last assertion is a consequence of Remark \ref{ Remark. 4.15}  and Definition \ref{Def. 4.13}. }
\end{Example}

\section{The  general case }

 Now, we go back to the general case of an ideal not necessarily
 generated by the maximal minors of an $m\times(m+1)$-matrix.

Let $f \in \I_d$ be any form that can be included in a standard
basis  \basis, or, equivalently, that does not lie in $\I_{d-1}S$.
As we just noticed, the fact that $f$ is essential depends on the
basis \basis. Our aim is to investigate how the {\it  nature of f
with respect to essentiality} (briefly: {\it the nature of f})
changes with  \basis.
Some lemmas will be useful.
\medskip
\begin{Lemma} \label{Lemma. 2.1} {\rm Let ${\cal B} = ( f_1,...,f_{i-1}, f_i, f_{i+1},...,f_m)$, ${\cal B'} = ( f'_1,...,f'_{i-1}, f_i+h, f'_{i+1},...,f'_m)$ be two standard bases of $\I$, where $h,f'_j \in \H=(f_1,...,\check{f_i},...,f_m)$. Then $f= f_i \in \I_d$ is essential (resp. inessential) with respect to ${\cal B}$ iff $f+h$ is so with respect to ${\cal B}'$.}
\end{Lemma}

\medskip \noindent
{\bf Proof}

It is enough to observe that: $ \H_{(f,\ {\cal B})} = \H_{(f+h,\ {\cal B}')}$. As a consequence:

$$(f+h) \M^t \subset \H_{(f+h,\ {\cal B}' )}\Leftrightarrow f\M^t \subset \H_{(f,\ {\cal B})}.$$

\finedim

\medskip

 Hence a basis change acting only on
\basis$- \{ f\}$ does not modify $f$'s nature, as it does not
modify $ \H = ({\cal B}(\I) -\{f\})S$.  In particular, that
happens for a basis change acting on elements of degree different
from $d$.

\medskip

\begin {Lemma} \label{Lemma. 2.2} {\rm The nature of $f=f_i \in {\cal B}(\I) = (f_1,...,f_m)$,  with respect to another basis ${\cal B}'(\I)$ containing it, is the same it has with respect to a basis of the type $ {\cal \widetilde{ B}}= (f_j + a_j f), j=1...m,  a_i = 0$, with the $a_j$ properly chosen.}
\end{Lemma}

\medskip\noindent {\bf Proof}

Let ${\cal B}'(\I)= (f'_1,..., f'_{i-1}, f, f'_{i+1},...,f'_m)$ be any other basis, linked to \basis \   by the relation ${\cal B}'(\I)={\cal B}(\I) T$, where $T$ is an invertible matrix. The submatrix $T_{ii}$, obtained from $T$ by deleting the row and the column of index $i$, is still invertible and acts on  ${\cal B}_1 = (f_1,...,\check{f}_i,...,f_m$. The base change acting on ${\cal B'(\I)}$ with $T_{ii}^{-1}$  produces a basis $\widetilde{\cal{B}}$ as described in the statement and the nature of $f$ with respect to it is the same that it had with respect to ${\cal B}'(\I)$, thanks to Lemma \ref{Lemma. 2.1}.
\finedim

\medskip

The previous Lemma gives immediately the following

\begin{Proposition} \label{Prop. 2.13'} {\rm To decide the nature of $f =f_i \in {\cal B} = (f_1,...,f_m)$, when ${\cal B}$ is replaced by any ${\cal B}'$ containing it, it is enough to consider just the bases ${\cal \widetilde{ B}}$ obtained from  ${\ B}$ by replacing $f_j$ with $ (f_j + a_{j}f), j \not=i$, for all (degree-allowed) forms $a_j$.}
\end{Proposition}
\medskip
{\bf Remark} \quad Let us observe that Corollary \ref{Cor. 2.11} can be viewed as a consequence of the previous Proposition, as the replacement of $C_r$ with $C'_r$ corresponds to a replacement of $g_i$ with $g_i-t_{ri}g_r, i=1...m,i\not=r .$

\bigskip
\begin{Lemma}\label{Lemma 2.28} {\rm  Let $c, c_1$ be elements of \basis \ in the same degree $d$, both inessential (resp. essential), and \ $c$ \ s.i. (resp. s.e.). Then a replacement of $c$ with $c+\alpha c_1$ cannot change the nature of $c_1$.}
\end{Lemma}
\medskip \noindent
{\bf Proof}

Let us consider the new basis ${\cal B}'(\I)$ in which $c$ is replaced by $c+\alpha c_1=c', \alpha \not= 0$. The nature of its element $c_1$ does not change if we replace it with $c_1-\alpha^{-1} c' = -\alpha^{-1} c$; this shows that $c_1$ preserves its former nature with respect to ${\cal B}'(\I)$.
\finedim

\bigskip
Now we point our attention on the bases having the greatest
(respectively: smallest) number of essential generators in a
chosen degree $d$: let us denote ${\cal B}^{(d)}(\I)$ ( resp.
${\cal C}^{(d)}(\I)$) any of them and $\nu_e(d)$ (resp. $\mu_e(d)$)
the number of their essential entries in degree $d$.

\begin{Proposition} \label{Prop.  2.14} {\rm There exist bases
${\cal B}_{Max}(\I)$ (resp. ${\cal B}_{min}(\I)$ ) having, in
every degree $d$, exactly $\nu_e(d)$ (resp. $\mu_e(d)$) essential
generators.}

\end{Proposition}

\medskip\noindent {\bf Proof}

We will prove the statement for ${\cal B}_{Max}$ ; the same
reasoning can be repeated for ${\cal B}_{min}$. As usual,
$\alpha$ is the minimal degree of an element of $\I$. Let us
denote ${\cal B}_M^d$ a basis satisfying the required condition
for every degree $\leq d$. We will prove the existence of a ${\cal
B}_M^d$, for every $d$, using induction on $d$. For $d=\alpha$, we
can chose ${\cal B}_M^{\alpha} = {\cal B}^{(\alpha)}$, for some choice of ${\cal B}^{(\alpha)}$. Now, let us
suppose the existence of a ${\cal B}_M^d$ and produce a ${\cal
B}_M^{d+1}$. To obtain ${\cal B}_M^{d+1}$ it is sufficient to
replace in ${\cal B}^{(d+1)}$ the part of degree $\leq d$ with the
analogous of the chosen ${\cal B}_M^d$ and modify the generators of larger degree as follows.
 Let us denote $\phi_i$ any element of $B^d_M$ of degree $>d$ and $\psi_i$ any element of ${\cal B}^{(d+1)}$ of degree $>d$. We can write $ \psi_i = \sum _j a_j \phi_j + \delta_i$, where $\delta_i \in \I_d S$. Let us set $ \psi'_i = \psi_i - \delta_i = \sum_j a_j \phi_j$.  We claim that we obtain a ${\cal B}_M ^{d+1}$  by replacing the generators of $B^d_M$ of degree $\geq d+1$ with the $\psi'_i$. In fact, with respect to this basis, any $\psi'_i$ has the same nature of the corresponding $\psi_i$ with respect to ${\cal B}^{d+1}$, so that in degree $d+1$ we have the maximum number of essential elements; moreover, the elements of degree $\leq d$ have the same nature with respect to ${\cal B}_M^d$ and with respect to ${\cal B}_M^{d+1}$, as the change we made in degree $>d$ does not involve them.
 \finedim

\medskip
\begin{Definition} \label {Def. 2.15} {\rm Every basis satisfying
the condition of Proposition \ref{Prop. 2.14} will be called {\it
maximal}  (resp. {\it minimal }) with respect to essentiality  or,
briefly, {\it e-maximal basis} (resp. {\it e-minimal basis}). Its
number of essential elements will be denoted $\nu_e(\I)$  (resp
$\mu_e(\I)$.}

\end{Definition}

\medskip
As a consequence of Proposition \ref{Prop. 2.9} we can state the
following

\begin{Corollary} \label{Cor. 2.16} {\rm If $\I \subset K[y_0,...,y_n], n
\geq 2$ is a perfect height $2$ ideal satisfying the condition:
$\nu(\I) \leq n+1$ , then every \basis \ is an e-maximal basis (more precisely, no basis contains inessential elements).}
\end{Corollary}

\noindent\medskip
 {\bf Proof} \quad As every Hilbert matrix of $\I$ has at most $n$ rows, it is enough to observe that $\M^h$
cannot be contained in an ideal generated by at most $n$ forms.
\finedim
\medskip \noindent
{\bf Remark} \quad Thanks to Dubreil's Theorem
saying that $\nu(\I) \leq \alpha(\I)+1$,
the condition of Corollary \ref{Cor. 2.16} is necessarily verified
if $\alpha(\I) \leq n$, that is if the minimal degree of a
hypersurface containing the corresponding scheme is $\leq n$.

\bigskip

Now we come back to the general situation.

\begin{Proposition} \label{Prop. 4.16} {\rm Let \basis \ and ${\cal B}'(\I)$ be two basis such that, in degree $d$, all their inessential  (resp. essential) elements are s.i. (resp.  s.e). Then, the subspace of $\I_d/(\I_{d-1}S_1)$ generated by the inessential (resp. essential) elements of \basis \ coincides with the one generated by the inessential (resp. essential) elements of ${\cal B}'(\I)$. }

\end{Proposition}

\noindent\medskip
 {\bf Proof}

 Let us consider the elements of \basis \ and ${\cal B}'(\I)$ in degree $d$:

 $${\cal B}_d(\I) = ( b_1,...,b_h, c_1,...,c_k);
\qquad
 {\cal B}'_d(\I) = ( b'_1,...,b'_h, c'_1,...,c'_k),$$

 where $b_i,b'_i$ are essential and $c_j,c'_j$ strongly inessential.

 In $\I/\I_{d-1}$ we have:  $$c'_j =\sum_{i=1}^k \gamma_i c_i  + \sum_{i=1}^h \beta_i b_i.$$

Let us prove that $\beta_i=0, i=1,...,h$. As we can exchange the role of \basis \ and ${\cal B}'(\I)$, that will be enough to complete the proof. So, let us suppose $\beta_i \neq 0$ for some $i$ and get a contradiction. In fact $\beta_i \neq 0$ implies that, in \basis, \ $b_i$ can be replaced by $c'_j$ without changing its nature, against the hypothesis that $c'_j$ is inessential with respect to every basis containing it.

Interchanging {\it inessential} and {\it essential} we prove the other part of the statement.
\finedim

Proposition \ref{Prop. 4.16} gives immediately the following consequences:

\begin {Corollary} \label{Cor. 4.17} {\rm Two basis whose inessential (resp. essential) elements are $s.i.$ (resp. $s.e.$) must have the same number of inessential (resp. essential) elements, degree by degree.}

\end {Corollary}

\begin {Theorem} \label {Th. 4.18} {\rm A standard basis is e-maximal (resp. e-minimal) iff its inessential (resp. essential) elements are strongly inessential (resp. strongly essential).       }

\end{Theorem}

\noindent\medskip
 {\bf Proof}

 We prove the statement for the e-maximal case, as the e-minimal one is analogous.

 First we prove that the inessential elements of an e-maximal basis ${\cal B}_M$ are s.i.. Let $c\in
{\cal B}_M$ be inessential; thanks to Proposition \ref{Prop. 2.13'}, it is enough to check that $c$ is still inessential with respect to the basis obtained from ${\cal B}_M$ by replacing each of its elements different from $c$, say $f_i$, with $f_i+a_i c= f'_i$. After such a replacement the nature of $f'_i$ is the same as the one of $f_i$, so that a change of nature of $c$ would imply the existence of a basis with $h+1$ essential elements, against the maximality of ${\cal B}_M$.

Viceversa, let  ${\cal B}$ be any basis whose inessential elements are s.i.. We just proved that every e-maximal basis ${\cal B}_M$ has such a property, so that Corollary \ref{Cor. 4.17} states that ${\cal B}$  and ${\cal B}_M$ have the same number of s.i. elements. As a consequence, also ${\cal B}$ is e-maximal.
\finedim

\medskip

Now, let us give a construction of an e-maximal (resp. e-minimal) basis.

\begin{Proposition}\label{Prop. 4.18} {\rm Starting from any basis \basis \ it is possible to produce an e-maximal basis (resp. e-minimal basis) containing all the s.i. (resp. s.e.) elements of \basis.}
\end{Proposition}

\noindent\medskip
 {\bf Proof}

Let us consider first the case of an e-maximal basis.

Thanks to Theorem \ref{Th. 4.18}, the aim is to produce a basis whose inessential elements are s.i.. So, we start to consider the inessential generators, non s.i., of lowest degree,
 following the order in which they appear in \basis: \ let $c_1$, deg $c_1=d$, be the first of them.  The replacement of some other elements $f_i\in {\cal B}(\I)$ with $f'_i = f_i+a_ic_1$  makes $c_1$ essential ( Lemma \ref{Lemma. 2.2}),  while $f'_i$, with respect to the new basis, has the same nature of $f_i$ (Lemma \ref{Lemma. 2.1}). We observe that the s.i. elements of degree $d$ are not involved, thanks to Lemma \ref{Lemma 2.28}. Let us denote ${\cal B}^1(\I)$ the new basis at this step, in which $c_1$ is  essential and, as a consequence, the number of inessential, but non s.i., elements is decreased . Then we go on dealing with ${\cal B}^1(\I)$ just as we did with ${\cal B}(\I)$. After a finite number of steps, we get a basis ${\cal B}^u(\I)$ whose inessential elements are s.i..

Analogously, it is possible to produce an e-minimal basis, starting from any basis  \basis: it is enough to replace $inessential$ with $essential$ in the previous construction.

\finedim

\bigskip

Now we turn our attention to the dehomogenization $\I_*$ of $\I$ with respect to a generic linear form $L$ and to the system of generators ${\cal B}_*$ obtained from \basis \ dehomogenizing every form appearing in it. In general, ${\cal B}_*$ is not a basis and our aim is to find its subsets that are bases and, among them, the ones of minimal cardinality.

The definition of inessential element can be generalized as follows.

\begin {Definition} \label { Def.4.1} {\rm A subset $T$ of \basis  \ is inessential iff $\I^{sat} = ({\cal B}(\I) - T)^{sat}$.}
\end {Definition}

\smallskip

Let us observe that if $T = \{t \}$, then $T$ is an inessential subset of \basis \ iff $t$ is inessential as an element.

\smallskip

Given a standard basis \basis, \ we point our attention on its maximal (with respect to $\subset$ ) inessential sets. Their interest lies on the following statement.

\begin{Proposition}\label{Prop.4.19} {\rm Let ${\cal B}(\I) =(b_1,...,b_h,c_1,...,c_k)$, where $T=(c_1,...,c_k)$ is a maximal inessential set. Then, in any dehomogenization $\I_*$ of $\I$ with respect to a linear form $L$, the set ${\cal B}_L(\I_*)= (b_{1*},...,b_{h*})$  is a set of generators of $\I_*$. Moreover, in the $K$-space $\M_1$, the subset of the linear forms $L$ such that ${\cal B}_L(\I_*)$ is not a basis of $\I_*$, is a finite union of proper linear subspaces (briefly, we can say that ${\cal B}_L(\I_*)$ is "generically" a basis).}
\end{Proposition}

\noindent\medskip
 {\bf Proof}

The first assertion comes immediately from the definition of inessential set. In fact, $c_i\in (b_1,...,b_h)^{sat}$ means that:
$$\forall L \in \M_1, (\exists t) c_i L^t \in(b_1,...,b_h).$$
As a consequence, $c_{i*} \in (b_{1*},...,b_{h*})$, in the dehomogenization with respect to $L$.

The second part of the statement can be proved just observing that \ $ b_{i^*} \in (b_{1*},...,b_{h*})$ \ is equivalent to \ $b_i L^t \in (b_1,...\check{b}_i,...,b_h)$, for some $t$. As
$$b_i L^{t_j}_j \in  (b_1,...\check{b}_i,...,b_h), j=1,2  \Rightarrow b_i(L_1+L_2)^{2 sup(t_1,t_2)} \in (b_1,...\check{b}_i,...,b_h),$$
the set of all the linear forms $L$ for which ${\cal B}_{*L}$ is not a basis is a finite union of linear subspaces $V_i$ of $\M_1$. The equality $V_i = \M_1$ cannot hold, as it should imply $b_i \in (b_1,...\check{b}_i,...,b_h)^{sat}$, against the hypothesis on the maximality of $T$.
\finedim

 The following proposition, lying on Proposition \ref{Prop. 2.17}, gives conditions equivalent to the one defining an inessential set $T = \{c_1,...,c_k\}$.

\begin{Proposition} \label{Prop. 4.14} {\rm Let \basis = $(b_1,...,b_h,c_1,...,c_k),$ $c_i$ inessential, ${\rm deg} c_i \leq  {\rm deg} c_{i+1}, i=1,...,k $. The following facts are equivalent:

$i)$ For any $ i=1,...,k, \ c_i$ is inessential with respect to the basis ${\cal B}_i = (b_1,...,b_h,c_1,...,c_i)$, generating an ideal $\A_i$.

$ii)$ If $\A$ is the ideal generated by $(b_1,...,b_h)$, then $\I^{sat}=\A^{sat}$.

$iii)$ For any $i=1,...,k$, for every form  \  $\alpha_{ij}, \ {\rm deg}\alpha_{ij} = {\rm deg} c_j- {\rm deg} c_i, \ j=i+1,...,k, \ c_i$ is inessential with respect to  ${\cal B}(\alpha_{ij})  = (b_1,...b_h,c_1,...,c_i,c_{i+1}+\alpha_{i,i+1}c_i, ...,c_k+\alpha_{ik}c_i)$.

$iv)$ For any $i=1...k, \ c_i$ is inessential with respect to every standard basis of the type
$(b_1,...,b_h,c_1,...,c_i, f_{i+1},...,f_k)$.}

\end{Proposition}

\noindent\medskip
 {\bf Proof}

$i)\rightarrow ii)$ It is enough to observe that $c_{i+1}\in \A_i^{sat}$ implies $(\A_{i+1})^{sat} = \A_i^{sat}$.

$ii)\rightarrow i)$ The hypothesis implies that $c_i$ is inessential with respect to $(b_1,...,b_h,c_i)$ and, as a consequence, with respect to $(b_1,...,b_h,c_1,...,c_i)$.

$i)\rightarrow iii)$ Obvious, because  $(b_1,...,b_h, c_1,...,c_{i-1)}S \subset (b_1,...,b_h,c_1,...,c_{i-1},c_{i+1}+\alpha_{i+1},c_i,...,c_k+\alpha_{ki}c_i)S$.

$ iii) \rightarrow i) $ We use induction on $k-i$.

If $k-i = 0$ both conditions say that $c_k$ is inessential with respect to \basis.

So, let us suppose the implication true until $k-(i+1)$ and prove it for $k-i$. Induction says that $c_{i+1},...,c_k \in (b_1,...,b_h,c_1,...,c_i)^{sat}$. So, it is enough to prove that $c_i\in((b_1,...,b_h,c_1,...,c_{i-1})^{sat}$ or, equivalently, that $c_i$ is inessential with respect to ${\cal B}_i$. If not, according to Proposition \ref{Prop. 2.17} the essentiality of $c_i$ with respect to ${\cal B}_i$ would mean that there exists a set $\{L_1,...,L_n, N\}$ of linear forms, generating $\M$, such that $N$ is regular both for $S/(\A_{i-1})^{sat}$ and for $S/(\A_i)^{sat}$ and moreover:

 \begin{equation}
 \forall t\in \N,  c_i N^t \notin ({\cal B}_{i-1}, c_i L_1,...,c_i L_n), t\gg 0, \label{eq. 4.1}
\end{equation}
Hence, we will prove that it is possible to find $v_j\in K, j=i+1,...,k$, such that

\begin{equation}
   c_i N^t \notin ({\cal B}_{i-1},c_{i+1}+ v_{i+1}N^{u_{i+1}}c_i,... , c_k+ v_k N^{u_k}c_i, c_i L_1,...,c_i L_n), t\gg 0, \label{eq. 4.2}
\end{equation}
so that $c_i$ is essential with respect to ${\cal B}(\alpha_{ij})$, where $\alpha_{ij}=v_jN^{u_j}c_i$, against condition $iii)$.

To this aim, it is enough to prove that the $v_j$'s can be chosen to realize the inclusion:

$((c_j+v_jN^{u_j} c_i)S)_t \subseteq ({\cal B}_{i-1}, c_i L_1,...,c_i L_n), t\gg 0, j>i,$

or equivalently:

\begin{equation}
  (c_j+v_jN^{u_j} c_i)L_w^t \in ({\cal B}_{i-1}, c_i L_1,...,c_i L_n), t\gg 0, j>i, w=1,...,n. \label{eq. 4.3}
\end{equation}

\begin{equation}
  (c_j+v_jN^{u_j} c_i)N^t \in({\cal B}_{i-1}, c_i L_1,...,c_i L_n), t\gg 0, j>i.\label{eq. 4.4}
\end{equation}

Now,(\ref{eq. 4.3}) is equivalent to : $ c_jL_w^t \in({\cal B}_i, c_i L_1,...,L_n), t\gg 0, w=1,...,n$. But we already observed that induction implies $c_j \in (b_1,...,b_h,c_1,...,c_i)^{sat}, j>i$, hence (\ref{eq. 4.3}) is true.

Let us consider (\ref{eq. 4.4}). Using induction, we have: $\exists t_0, c_j N^{t_0} = h_1 + c_iP + \alpha c_i N^{t_1}$, where $h_1 \in \A_{i-1}, P\in(L_1,...,L_n), \alpha \in K$. So, to realize (\ref{eq. 4.4}), it is enough to choose $v_j=-\alpha$.

$iii)\leftrightarrow iv)$ It is enough to use Lemma \ref{Lemma. 2.2}.

\finedim
\medskip

An immediate consequence of the previous Proposition  and of Theorem \ref {Th. 4.18} is the following:

\begin{Corollary} \label{Cor.4.15} {\rm  If \basis \ is an e-maximal basis, then the set of all its inessential elements is an inessential set.}
\end{Corollary}

\bigskip

As a consequence of Proposition \ref{Prop.4.19}, we can say that a subset of a standard basis \basis \ gives rise to a minimal basis, for the dehomogenization of $\I$ with respect to a generic linear form, iff it is the complement, in \basis,\ of a maximal inessential subset with maximal cardinality.

The following proposition says that we can find a maximal inessential set of maximal cardinality among the inessential sets of the e-minimal bases.

\begin{Lemma}{\label{Lemma 5.3} }{\rm Let \basis= $(b_1,...,b_h,c_1,...,c_k)$, where $T=(c_1,...,c_k)$ is an inessential set. In a basis: ${\cal B'}(\I) = (b'_1,...,b'_h,c'_1,...,c'_k)$, where  $(b'_1,...,b'_h)S = (b_1,...,b_h)S = \J$, the subset $T'=(c'_1,...,c'_k)$ is still inessential.}\end{Lemma}.

\noindent\medskip
 {\bf Proof}

 The hypothesis says that $T \subset \J^{sat}$. As a consequence, also the inclusion $T'\subset \J^{sat}$ holds.

 \finedim

 \begin{Proposition}\label{Prop.5.7}{\rm Let ${\cal B}_m(\I)$ be the e-minimal basis produced for \basis \ according to Proposition \ref{Prop. 4.18}. For every inessential subset $V \subset {\cal B}(\I)$, there exists an inessential subset $V' \subset {\cal B}_m(\I)$ with the same cardinality of $V$.}\end {Proposition}

 \noindent\medskip
 {\bf Proof}

 According to Proposition \ref{Prop. 4.18}, passing from \basis \ to ${\cal B}_m(\I)$, every element $f\in{\cal B}(\I)$ is replaced by $f'=f+\sum a_jf_j$, where the $f_j$'s are elements outside $V$. As a consequence, $V$ is replaced by $V'$, with its same cardinality, and we can apply Lemma {\ref {Lemma 5.3}, where $T=V$.
\finedim

From the previous Proposition it is immediate to get:

\begin {Corollary}\label{Cor.5.8} {\rm The maximal cardinality of the inessential subsets of ${\cal B}_m(\I)$ is not less than the one of the inessential subsets of \basis.}
\end{Corollary}

\bigskip

In the special case of perfect height $2$ ideals, Proposition \ref{Prop. 2.9}, Corollary
\ref{Cor. 2.11} and Proposition \ref{Prop. 4.14}  give rise to the following:

\begin{Proposition}\label{Prop.5.8} {\rm Let $\M(\I)$ be a Hilbert matrix with respect to the basis \basis. \ The subset $T$, $|T|=s$, of \basis\, corresponding to the columns $C_{i_1},...,C_{i_j},..., C_{i_s},\ i_1<i_2<...<i_s$, of $\M(\I)$, is inessential iff, $ \forall \ i_j \in (i_1,...,i_s)$, the following condition is satisfied:

(*) For every choice of the forms $t_{i_j,i_h}, h\geq j, \ deg \ t_{i_j,i_h} = deg \ g_h - deg\ g_j$, the entries of
$$C'_{i_j} = \sum_{h\geq j}  t_{i_j,i_h} C_{i_h},\quad t_{i_j,i_j}= 1$$
generate an ideal $\I_{C'_{i_j}}$ containing some power of the irrelevant ideal $\M$.}
\end{Proposition}

\bigskip

Now we point our attention on the bases whose inessential elements form an inessential set.

\begin
{Definition} \label {Def.4.2} {\rm A basis \basis \ whose essential elements generate an ideal $\E$ such that $\E^{sat} = \I^{sat}$  is called {\it essential basis} (briefly: {\it E-basis}).}
\end{Definition}

Now, our aim is to produce, starting from any standard basis \basis, an E-basis ${\cal B}_E(\I)$, containing all the essential elements of \basis.

\begin{Lemma}\label{Lemma 5.9} {\rm Let \basis $=(b_1,...,b_h,c_1,...,c_k)$, $b_i$ essential, $i=1...h$ , $c_i$ inessential, $i=1,...,k$ , be any standard basis. There exists a maximal inessential set $V \subset {\cal B}(\I)$ such that:

\begin{equation}
c_i \notin V \Rightarrow \{ c_i, V\cap (c_{i+1},...,c_k)\}
{\rm \ is \ not \ inessential.} \label {eq.5.1}
\end{equation}}
\end{Lemma}

\noindent\medskip
 {\bf Proof}

 We define $V$, step by step, by means of the following conditions:

 i) $c_k \in V$,

 ii) $c_i \in V, \ i<k \Leftrightarrow c_i \in ( {\cal B(\I)} - \{ c_i, V\cap (c_{i+1},...,c_k)\})^{sat}$.

 Condition  (\ref{eq.5.1}) coincides with $ii)$ and the inessentiality of $V$ is an immediate consequence of the definition of inessential set. Moreover, condition (\ref{eq.5.1}) implies the maximality of $V$.

\finedim

\begin{Proposition}\label{Prop.5.10} {\rm Starting from any basis \basis, it is possible to produce an E-basis ${\cal B}_E(\I)$, whose set of essential elements includes the ones of \basis.}
\end{Proposition}

\noindent\medskip
 {\bf Proof}

Let us point our attention on the maximal inessential set \ $V=(v_1,...,v_r)$ \ defined in Lemma \ref{Lemma 5.9}. Condition (\ref{eq.5.1}) says that every inessential element $c_i$ outside $V$ becomes essential by a replacement of $v_j$ with $v'_j = v_j + \sum a_{ij}c_i$, where the $a_{ij}$ are properly chosen step by step (see Proposition \ref{Prop. 4.14}). Moreover, Lemma \ref{Lemma 5.3} assures that the subset $V'=(v'_j), \ j=1...r$  is inessential in ${\cal B}'(\I)$, that, as a consequence, turns out to be an E-basis.
\finedim

\end{document}